# CENTRAL LIMIT THEOREM FOR A MANY-SERVER QUEUE WITH RANDOM SERVICE RATES

By Rami Atar

*Technion—Israel Institute of Technology*

Given a random variable $N$ with values in $\mathbb{N}$, and $N$ i.i.d. positive random variables $\{\mu_k\}$, we consider a queue with renewal arrivals and $N$ exponential servers, where server $k$ serves at rate $\mu_k$, under two work conserving routing schemes. In the first, the service rates $\{\mu_k\}$ need not be known to the router, and each customer to arrive at a time when some servers are idle is routed to the server that has been idle for the longest time (or otherwise it is queued). In the second, the service rates are known to the router, and a customer that arrives to find idle servers is routed to the one whose service rate is greatest. In the many-server heavy traffic regime of Halfin and Whitt, the process that represents the number of customers in the system is shown to converge to a one-dimensional diffusion with a random drift coefficient, where the law of the drift depends on the routing scheme. A related result is also provided for nonrandom environments.

**1. Introduction.** Many-server queues in heavy traffic have been the subject of much research, both for their theoretical interest and for their importance in practical applications (a review of some of the more theoretical results appears in Whitt [11]). Central limit theorem (CLT) results for such models are known under various settings, including ones where the servers are homogeneous and ones where they are heterogeneous, but models in which the service rates are random have not so far been treated in a many-server CLT framework. Incorporating random service rates into a many-server setting appears to be a natural way to model uncertainty in the abilities of individual servers. This approach is taken here in analyzing a many-server, single class model under two routing policies.

A CLT regime for a many-server queue was first studied in detail by Halfin and Whitt [5]. In this paper one considers a parameterization in which the number of servers and the rate of arrivals are scaled up at a (nearly) fixed proportion under which the queue remains critically loaded. It is assumed that the arrivals occur according to a renewal process and that the servers all have exponential service time distribution of the same rate, and it is shown that the second-order asymptotics of the process representing the number of customers in the system, is given as a one-dimensional diffusion. It is well understood that one-dimensional diffusions are not to be expected when one relaxes the exponential assumption. Indeed, Puhalskii and Reiman [8] prove that an $(n-1)$-dimensional diffusion appears in the limit when the service time distribution is of phase type with $n$ phases; Reed [9] and Kaspi and Ramanan [6] derive recursions in function space for general service time distribution (for both fluid and diffusion limits). Similarly, when one keeps the exponential assumption but relaxes the server homogeneity assumption—a situation which occurs, in particular, when service rates are random—one does not in general obtain a one-dimensional diffusion limit (we do not prove this claim). Since one-dimensional diffusions constitute a significant simplification that the model undergoes in the limit, it is desirable to understand when they occur. The focus of this paper is on models that feature exponential servers with random rates, and give rise to one-dimensional diffusions.

In the more conventional heavy traffic regime, where the arrival and service rates are scaled up but the number of servers is not, fluid and diffusion limit results have been proved in Choudhury et al. [4] for a random, time-varying environment setting. Besides the fact that the time-varying environment aspect is central in [4], the nature of these results is quite different from those of the current paper also in that in a many-server setting, the environment stochasticity that appears in the limiting diffusion process (in the form of a random drift coefficient) has an ingredient that originates from a CLT for the random service rates themselves. We also mention the fluid scale treatment of a many-server model in a setting of random environment, Whitt [10], where randomness enters in the number of servers and arrival rates.

The model under study has a random number $N$ of exponential servers with i.i.d. service rates $\mu_k$, $k = 1, \ldots, N$. When customers arrive into the system they are either queued in a buffer with infinite room, or routed to a server according to a specified routing policy. Customers from the queue are routed to servers according to a first-come-first-served rule. Each customer leaves the system when its service requirement is fully processed. The routing policies are *work conserving*, in the sense that no server may be idle when at least one customer is in the buffer. The service policy is *noninterruptible* (a term sometimes referred to as *nonpreemptive*), in the sense that each customer is assigned exactly one server, that then continuously processes its service requirement to completion. In the first routing policy under consideration, referred to as P1, the service rates $\{\mu_k\}$ need not be known to the router, and each customer to arrive at a time when some servers are idle is routed to the server that has been idle for the longest time, or otherwise it is queued. Under the second routing policy, P2, the service rates are known to the router, and a customer that arrives to find idle servers is routed to the one whose service rate is greatest. In our main results (Theorems 2.1 and 2.2) we find second-order asymptotics for the process that represents the number of customers in the system in the form of a one-dimensional diffusion with a random drift coefficient (that depends on the policy).

By its very definition, policy P1 expresses a form of fairness, because when several servers are free, the one selected for the next incoming job is the one that has been idle for the longest time. In addition, we will show that an asymptotic property of load balancing takes place, in the sense that when the load on the servers is relatively low so that some servers have to idle, they all experience idle periods of approximately the same length (as expressed by Proposition 2.1, which also identifies the asymptotic length of idle periods in terms of the limiting diffusion alluded to above).

Armony [1] analyzes policy P2 (referred to in [1] as Fastest Server First) in a deterministic environment, where the service rates take a finite number of values. Our result on policy P2 (Theorem 2.2) is analogous to Proposition 4.2 of [1], but a significant difference is the presence of a random drift term in the diffusion limit of our Theorem 2.2, originating from environment stochasticity.

The result regarding policy P1 (Theorem 2.1) may also be interesting in a deterministic environment setting, and thus we formulate such a version of it in Corollary 4.1.

We point out that the phenomena which make the diffusion limit one-dimensional are very different under policies P1 and P2. A heuristic discussion at the end of Section 2 explains the form of the limit processes in both cases.

The organization of the paper is as follows. Section 2 contains a description of the model, statement of the main results, namely Theorem 2.1 and Proposition 2.1 for the routing policy P1 and Theorem 2.2 for P2, and a heuristic discussion. The proofs appear in Section 3. Finally, in Section 4 we analyze policy P1 in deterministic environment.

**2. Model description and main results.** We fix some notation. For a positive integer $d$, we denote by $\mathbb{D}(\mathbb{R}^d)$ the space of functions from $\mathbb{R}_+$ to $\mathbb{R}^d$ that are right continuous on $\mathbb{R}_+$ and have finite left limits on $(0, \infty)$ (RCLL), endowed with the usual Skorohod topology [3]. For $X \in \mathbb{D}(\mathbb{R}^d)$ and $t > 0$, we denote $\Delta X(t) = X(t) - X(t-)$. If $X^n$, $n \in \mathbb{N}$ and $X$ are processes with sample paths in $\mathbb{D}(\mathbb{R}^d)$ (resp., real-valued random variables) we write $X^n \Rightarrow X$ to denote weak convergence of the measures induced by $X^n$ on $\mathbb{D}(\mathbb{R}^d)$ (resp., on $\mathbb{R}$) to the measure induced by $X$, as $n \to \infty$. For $X \in \mathbb{D}(\mathbb{R})$ we write $|X|_{*,t} := \sup_{0 \leq s \leq t} |X(s)|$. For a collection $\mathcal{A}$ of random variables, $\sigma\{\mathcal{A}\}$ denotes the sigma-field generated by this collection.

We now rigorously describe the model: the servers and their rates, the initial configuration, the arrival process, and finally the routing policies. A complete probability space $(\Omega, \mathcal{F}, P)$ is given, supporting all random variables and stochastic processes defined below. Expectation w.r.t. $P$ is denoted by $E$. The model is parameterized by $n \in \mathbb{N}$. For each $n$, let $N^n$ be a random variable with values in $\mathbb{N}$, representing the number of servers. Let $(\widetilde{\mu}_k, \widehat{\mu}_k)$, $k \in \mathbb{N}$ be $\mathbb{R}^2$-valued i.i.d. random variables, and let

$$(2.1) \qquad \mu_k^n := \widetilde{\mu}_k + n^{-1/2}\widehat{\mu}_k, \qquad k = 1, \ldots, N^n$$

represent the service rate of server $k$. It is assumed that $\mu_k^n$ are nonnegative. The marginal distribution of $\widetilde{\mu}_1$ (resp., $\widehat{\mu}_1$) is denoted by $\widetilde{m}$ (resp., $\widehat{m}$). It is also assumed that the random variables $\widetilde{\mu}_k$ are nonnegative, and that

$$(2.2) \qquad \mu := \int_{[0,\infty)} x \, d\widetilde{m} \in (0, \infty).$$

The random variables $\widehat{\mu}_k$ are assumed to be bounded, namely $P(|\widehat{\mu}_1| \leq \overline{\mu}) = 1$, for some constant $\overline{\mu} < \infty$. Also, the number of servers is assumed to be bounded by $2n$, that is, $P(N^n \leq 2n) = 1$, and to be asymptotic to $n$, in the sense that $N^n/n \Rightarrow 1$.

To describe the initial configuration, let $Q_0^n$ be a $\mathbb{Z}_+$-valued random variable, representing the initial number of customers in the buffer. Let $B_{k,0}^n$, $k = 1, \ldots, N^n$ be $\{0,1\}$-valued random variables representing the initial state of each server as follows: $B_{k,0}^n = 1$ if server $k$ initially serves a customer. It is assumed that $Q_0^n > 0$ if and only if $B_{k,0}^n = 1$ for all $k = 1, \ldots, N^n$. We denote $X_0^n = Q_0^n + \sum_{k=1}^{N^n} B_{k,0}^n$. This random variable represents the total number of customers initially in the system. By assumption, we have the relation $Q_0^n = (X_0^n - N^n)^+$.

The "second-order asymptotics" of the random variables $X_0^n$ and $N^n$, defined below, are assumed to satisfy

$$(2.3) \qquad (\widehat{X}_0^n, \widehat{N}^n) := (n^{-1/2}(X_0^n - N^n), n^{-1/2}(N^n - n)) \Rightarrow (\xi_0, \nu),$$

where $(\xi_0, \nu)$ is an $(\mathbb{R}^2)$-valued random variable.

The arrivals are modeled as renewal processes with finite second moment for the interarrival time. To this end, we are given parameters $\lambda^n > 0$, $n \in \mathbb{N}$ satisfying $\lim_n \lambda^n/n = \lambda > 0$, and a sequence of strictly positive i.i.d. random variables $\{\check{U}(l), l \in \mathbb{N}\}$, with mean $E\check{U}(1) = 1$ and variance $C_{\check{U}}^2 = \mathrm{Var}(\check{U}(1)) \in [0, \infty)$. With $\sum_1^0 = 0$, the number of arrivals up to time $t$ for the $n$th system is given by

$$A^n(t) = \sup\left\{l \geq 0 : \sum_{i=1}^{l} \frac{\check{U}(i)}{\lambda^n} \leq t\right\}, \qquad t \geq 0.$$

The "heavy traffic" condition assumed throughout, that makes the system critically loaded, relates the arrival and service rates as follows:

$$(2.4) \qquad \lambda = \mu.$$

The arrival rates are moreover assumed to satisfy the second-order relation

$$(2.5) \qquad \lim_n n^{-1/2}(\lambda^n - n\lambda) = \widehat{\lambda},$$

for some $\widehat{\lambda} \in \mathbb{R}$.

For each $k = 1, \ldots, N^n$, we let $B_k^n$ be a stochastic process taking values in $\{0,1\}$, representing the status of server $k$. When $B_k^n(t) = 1$ we say that server $k$ is busy. For $k = 1, \ldots, N^n$, let $R_k^n$ (resp., $D_k^n$) be a $\mathbb{Z}_+$-valued with nondecreasing RCLL sample paths, representing the number of routings of customers to server $k$ within $[0, t]$ (resp., the number of jobs completed by server $k$ by time $t$). Thus

$$(2.6) \qquad B_k^n(t) = B_{k,0}^n + R_k^n(t) - D_k^n(t), \qquad k = 1, \ldots, N, t \geq 0.$$

We let $I_k^n(t) = 1 - B_k^n(t)$ for $k = 1, \ldots, N^n$, and $t \geq 0$. The assumptions on $R_k^n$ depend on the routing policy, and will be specified later. To describe the processes $D_k^n$, let $\{S_k, k \in \mathbb{N}\}$ be i.i.d. standard Poisson processes, each having right-continuous sample paths. The processes $D_k^n$ are assumed to satisfy

$$(2.7) \qquad D_k^n(t) = S_k(T_k^n(t)), \qquad k = 1, \ldots, N^n$$

where

$$(2.8) \qquad T_k^n(t) = \mu_k^n \int_0^t B_k^n(s) \, ds, \qquad k = 1, \ldots, N^n.$$

The stochastic primitives introduced thus far are

$$(\{\widetilde{\mu}_k, \widehat{\mu}_k\}_{k \in \mathbb{N}}, \{N^n, X_0^n\}_{n \in \mathbb{N}}, \{B_{k,0}^n\}_{k=1,\ldots,N^n}^{n \in \mathbb{N}}), \qquad \{S_k\}_{k \in \mathbb{N}}, \qquad \{A^n\}_{n \in \mathbb{N}}.$$

We have further assumptions regarding their joint law. First, it is assumed that the three random objects above are mutually independent. Moreover, $\{\widetilde{\mu}_k, \widehat{\mu}_k\}_{k \in \mathbb{N}}$ and $\{N^n, X_0^n\}_{n \in \mathbb{N}}$ are assumed to be mutually independent. Note, in particular, that it is not assumed that $\{B_{k,0}^n\}$ are independent of $\{\mu_k^n\}$.

Let $X^n$, $Q^n$ and $I^n$ be defined as

$$(2.9) \qquad X^n(t) = X_0^n + A^n(t) - \sum_{k=1}^{N^n} D_k^n(t),$$

$$(2.10) \qquad Q^n(t) = Q_0^n + A^n(t) - \sum_{k=1}^{N^n} R_k^n(t),$$

$$(2.11) \qquad I^n(t) = \sum_{k=1}^{N} I_k^n(t).$$

These processes represent the number of customers in the system, the number of customers in the buffer, and, respectively, the number of servers that are idle. It is assumed that the routing policy is work conserving, in the sense that

$$(2.12) \qquad Q^n(t) = (X^n(t) - N^n)^+, \qquad I^n(t) = (X^n(t) - N^n)^-.$$

Finally we come to the routing mechanism. We denote by $K^n(t)$ the set of servers that are "available" at time $t$. This includes those that are idle at $t-$, and those that have just finished a job at time $t$. More precisely, let $K_0^n = \{k : B_{k,0}^n = 0\}$. Then $K^n(0) = K_0^n$, and

$$(2.13) \qquad K^n(t) = \{k : I_k^n(t-) = 1\} \cup \{k : \Delta D_k^n(t) = 1\}, \qquad t > 0.$$

Let, for $t > 0$,

$$(2.14) \qquad \begin{aligned} RT^n(t) &= \{s \in (0,t] : K^n(s) \neq \varnothing\} \\ &\quad \cap \{s \in (0,t] : Q^n(s-) > 0 \text{ or } \Delta A^n(s) = 1\} \end{aligned}$$

represent the set of routing times, that is, the set of times at which customers are routed to servers up to $t$. To define the first policy we need the following notation. Let

$$(2.15) \quad H_k^n(t) = \inf\{s \in [0,t) : I_k^n(u) = 1 \text{ for all } u \in [s,t)\} \wedge t, \qquad t > 0.$$

$H_k^n(t)$ is equal to the last time prior to $t$ when server $k$ became idle if $I_k^n(t-) = 1$, and it is equal to $t$ otherwise. Also let

$$(2.16) \qquad H^n(t) = \begin{cases} \min\{H_k^n(t) : k \in K^n(t)\}, & K^n(t) \neq \varnothing, \\ t, & K^n(t) = \varnothing, \end{cases}$$

and

$$(2.17) \qquad h^n(t) = t - H^n(t).$$

The policy that we now define routes a customer to the server that has been idle for the longest time since it last served, or to one of the servers that has been idle since time zero, if such servers exist. Thus, with

$$(2.18) \qquad \kappa^n(t) = \min\{k \in K^n(t) : H_k^n(t) = H^n(t)\}$$

if $K^n(t)$ is nonempty, and $\kappa^n(t) = 0$ otherwise, the routing processes are given as

$$(2.19) \qquad R_k^n(t) = \sum_{s \in RT^n(t)} 1_{\{\kappa^n(s) = k\}}.$$

All the stochastic processes above are now well defined. More precisely, given $n$ and the data $\{\mu_k^n\}$, $N^n$, $X_0^n$, $\{B_{k,0}^n\}$, $\{S_k\}$ and $A^n$, there exists a unique process $\bar{\Sigma}_1^n = (X^n, Q^n, I^n, K^n, RT^n, H^n, \kappa^n, \{B_k^n, R_k^n, D_k^n, T_k^n, H_k^n\})$ such that equations (2.6)–(2.19) are satisfied for all $t \geq 0$, a.s. The sample paths of the process $\Sigma_1^n = (A^n, X^n, Q^n, I^n, \{B_k^n, D_k^n, R_k^n\}_{k=1,\ldots,N^n})$ are a.s. piecewise constant and right-continuous, and, a.s., each time $t$ when one of the $N^n + 1$ processes $(\{D_k^n\}, A^n)$ jumps, the size of the jump is one, and no one of the other $N^n$ processes jumps at $t$. These facts can be proved by induction on the times at which one of the $N^n + 1$ processes alluded to above jump, in a straightforward way, and the elementary proof is omitted. We refer to the process $\bar{\Sigma}_1^n$ as *policy* P1.

We let $\widehat{X}^n$ be a centered, renormalized version of the process $X^n$, defined by
$$\widehat{X}^n = n^{-1/2}(X^n - N^n).$$

THEOREM 2.1. *Assume $\int x^2 \, d\widetilde{m} < \infty$. Then, under policy P1, the processes $\widehat{X}^n$ converge weakly to the unique solution $\xi$ to the equation*

$$(2.20) \qquad \xi(t) = \xi_0 + \sigma w(t) + \beta t + \gamma \int_0^t \xi(s)^- \, ds, \qquad t \geq 0,$$

*where $\sigma^2 = \lambda C_{\check{U}}^2 + \mu \equiv \mu(C_{\check{U}}^2 + 1)$, $\beta = \widehat{\lambda} - \widehat{\mu} - \zeta - \mu\nu$, $\widehat{\mu} = \int x \, d\widetilde{m}$, $\zeta$ is a normal random variable with parameters $(0, \int (x-\mu)^2 \, d\widetilde{m})$, $\gamma = \int x^2 \, d\widetilde{m} / \int x \, d\widetilde{m}$, the process $w$ is a standard Brownian motion, and the three random objects $(\xi_0, \nu)$, $\zeta$ and $w$, are mutually independent.*

REMARK 2.1. It is a standard fact that there exists a strong solution $\xi$ to (2.20), adapted to the augmentation of the filtration $\sigma\{\xi_0, \beta\} \vee \sigma\{w_s : s \leq t\}$ [7], Theorem V.3.7, and that it is strongly unique (uniqueness, in fact, holds in the class of all continuous processes).

REMARK 2.2. Let $\widehat{I}^n = n^{-1/2} I^n$ and $\widehat{Q}^n = n^{-1/2} Q^n$. Then by (2.12), it follows from Theorem 2.1 that $(\widehat{X}^n, \widehat{I}^n, \widehat{Q}^n) \Rightarrow (\xi, \xi^-, \xi^+)$, and thus the result regards also the limit behavior of the queue length and idleness processes.

The length of the idle periods experienced by the servers can also be analyzed, and as shown by the result below, can be expressed in terms of the diffusion process $\xi$. The result also shows that a certain form of fairness is achieved. Roughly speaking, if at a certain moment $s$ there are relatively few customers [in the sense that the limiting process satisfies $\xi(s) < 0$], then the load is balanced among all servers that at or near $s$ have become available, in such a way that they all experience an idle period of approximately the same length.

To state the precise claim, we need the following notation. Given $0 < s < t$, $[s, t)$ is said to be an idle period for server $k$ if at $s$ it has finished serving a customer, and at $t$ it has started serving its next customer. $t - s$ is the length of that idle period. If the above occurs with $t = s$, we regard it as an idle period of length zero. For $0 < s < t$, let $\mathbf{K}_{s,t}^n$ denote the set of all servers that at some time $u \in [s, t]$ have finished an idle period. Let $\mathbf{I}_t^n(k)$ denote the length of the last idle period of server $k$, completed at or before time $t$ (left undefined if no such idle period exists), and

$$\bar{\mathbf{I}}_{s,t}^n = \sup\{\mathbf{I}_t^n(k) : k \in \mathbf{K}_{s,t}^n\}, \qquad \underline{\mathbf{I}}_{s,t}^n = \inf\{\mathbf{I}_t^n(k) : k \in \mathbf{K}_{s,t}^n\}.$$

PROPOSITION 2.1. *Under the assumptions of Theorem 2.1, given $s > 0$ and a sequence $\{t_n\}$ with $t_n > s$, $t_n \to s$ and $n^{1/2}(t_n - s) \to \infty$, one has*

$$n^{1/2}\bar{\mathbf{I}}_{s,t_n}^n \Rightarrow \mu^{-1}\xi(s)^-, \qquad n^{1/2}\underline{\mathbf{I}}_{s,t_n}^n \Rightarrow \mu^{-1}\xi(s)^-.$$

REMARK 2.3. One can think of the policy under study as if servers that become available enter a queue and are assigned new jobs according to a FIFO discipline. Viewed this way, the above result expresses a form of sample-path Little's law. Indeed, the number of idle servers is asymptotic (in diffusion scale) to $\xi^-$, the average rate at which servers become idle is (in fluid scale) $\mu$, and the time servers "wait for a job" is asymptotic (at diffusion scale) to the ratio $\mu^{-1}\xi^-$.

The second routing scheme that we consider is defined as follows. It is assumed that the servers are re-ordered in such a way that $\mu_1^n \leq \mu_2^n \leq \cdots \leq \mu_{N^n}^n$. Recall the set-valued process $K^n$, representing the set of available servers (2.13). We let

(2.21) $$\kappa^n(t) = \max\{k \in K^n(t)\},$$

if $K(t) \neq \varnothing$, and $\kappa(t) = 0$ otherwise [equation (2.21) is to be considered in place of (2.18)]. As before, let $RT^n$ and $R_k^n$ be defined via (2.14) and (2.19). The policy under consideration thus selects for each customer the fastest server available at the time of routing. The precise definition of the policy is via equations (2.6)–(2.14), (2.19) and (2.21), and the processes

$$\bar{\Sigma}_2^n = (X^n, Q^n, I^n, K^n, RT^n, \kappa^n, \{B_k^n, R_k^n, D_k^n, T_k^n\})$$

and

$$\Sigma_2^n = (A^n, X^n, Q^n, I^n, \{B_k^n, D_k^n, R_k^n\}_{k=1,\ldots,N^n}).$$

Assertions analogous to those stated after equation (2.19), regarding the processes $\bar{\Sigma}_1^n$, $\Sigma_1^n$ and equations (2.6)–(2.19), hold true for the processes $\bar{\Sigma}_2^n$, $\Sigma_2^n$ and equations (2.6)–(2.14), (2.19) and (2.21). The process $\bar{\Sigma}_2^n$ is referred to as *policy* P2.

For simplicity, we shall assume that the initial arrangement of customers is such that the faster servers are all busy, in the sense that

$$B_{k,0}^n = 1_{\{k > I_0^n\}}, \qquad (2.22)$$

where $I_0^n = (X_0^n - N^n)^-$ is the initial number of idle servers. Denote $\mu_{\min} = \operatorname{ess\,inf} \widetilde{m}$ and $\mu_{\max} = \operatorname{ess\,sup} \widetilde{m}$.

THEOREM 2.2. *Assume $\mu_{\max} < \infty$. Let the initial configuration satisfy (2.22). Then, under policy* P2, *the processes $\widehat{X}^n$ converge weakly to the unique solution $\xi$ to the equation*

$$\xi(t) = \xi_0 + \sigma w(t) + \beta t + \mu_{\min} \int_0^t \xi(s)^- \, ds, \qquad t \geq 0, \qquad (2.23)$$

*where $(\sigma, \beta, w)$ are as in Theorem 2.1.*

Note that Remarks 2.1 and 2.2 apply to policy P2 as well.

We end this section with a heuristic explanation of the form of the limiting processes. By calculation [cf. (3.8)–(3.12)], the process $\widehat{X}^n$ satisfies

$$\widehat{X}^n(t) = \widehat{X}_0^n + W^n(t) + b^n t + F^n(t), \qquad (2.24)$$

where

$$F^n(t) = n^{-1/2} \sum_{k=1}^N \mu_k \int_0^t I_k(s) \, ds.$$

For the policies under consideration, we show that $W^n$ converges to a Brownian motion and $b^n$ to a random variable. The term $F^n$ expresses decrease of processing rate due to idleness. It is the limit behavior of this term that differs from one policy to another, and whether it can be well approximated by a quantity of the form $\int_0^\cdot c(\widehat{X}^n(s)) \, ds$ will determine if the limit is a one-dimensional diffusion. To give an idea on why $F^n$ behaves as stated in the above results, consider first policy P1 and let $I(t; \Delta x)$ denote the number of servers with $[x, x + \Delta x)$-valued service rate, that are idle at time $t$, where $\Delta x$ is relatively small. Because the service time distribution of the relevant servers is exponential with rate $x$ approximately, and their number is roughly $N\widetilde{m}([x, x + \Delta x))$, the number of servers with $[x, x + \Delta x)$-valued service rate that become idle on a time interval $[a, b)$ is roughly Poisson with parameter $N\widetilde{m}([x, x + \Delta x))x(b - a)$. Therefore a first-order approximation for $I$ is

$$I(t; \Delta x) \approx N\widetilde{m}([x, x + \Delta x))x h^n(t). \qquad (2.25)$$

With $N \approx n$, we have

$$n^{-1/2} \sum_{k=1}^N \mu_k I_k(t) \approx n^{-1/2} \int x I(t; dx) \approx n^{1/2} h^n(t) \int x^2 \widetilde{m}(dx).$$

Similarly

$$\widehat{I}^n(t) = n^{-1/2}\sum_{k=1}^{N} I_k(t) \approx n^{1/2} h^n(t) \int x\widetilde{m}(dx)$$

hence

$$F^n(t) \approx \frac{\int x^2 \widetilde{m}(dx)}{\int x\widetilde{m}(dx)} \int_0^t \widehat{I}^n(s)\,ds = \gamma \int_0^t X^n(s)^-\,ds,$$

which, along with (2.24) explains the limit (2.20).

A rigorous version of the rough equality (2.25), on which this heuristic argument is built, can be found in the proof of Theorem 2.1 [i.e., (3.15) and the proven smallness of the (integrated) error terms in this equation].

For policy P2, the selection of free servers with higher rates encourages idleness among the servers with lowest rates, and this results in

(2.26) $$I(t;\Delta x) = o(n^{1/2}) \qquad \text{for any } x > \mu_{\min},$$

hence in the approximate equality

$$F^n \approx \mu_{\min} \int_0^{\cdot} \widehat{I}^n(s)\,ds = \mu_{\min} \int_0^{\cdot} \widehat{X}^n(s)^-\,ds.$$

This explains the limit (2.23). A rigorous statement of (2.26) appears in (3.26) as a part of the proof of Theorem 2.2. This phenomenon is similar to what appears in [1] in the case of a deterministic environment and a finite number of server pools, and the reader is referred to this reference for further discussion.

**3. Proofs.** For simplicity, we will omit the symbol $n$ from the notation of all random variables and stochastic processes. There will be no confusion with variables and processes that do not depend on $n$. The deterministic parameters that depend on $n$ will still have superscript $n$ in their notation.

The proposition below provides a tool by which one can replace the Poisson processes that drive the service processes of the individual servers by ones that jointly drive service processes of collections of servers.

PROPOSITION 3.1. *Fix $n \in \mathbb{N}$. Let $(K_1,\ldots,K_q)$ be a partition of $(1,\ldots,N)$ that is measurable on $\sigma\{N,\{\widetilde{\mu}_k,\widehat{\mu}_k\}\}$. Let $\{S^{(1)},\ldots,S^{(q)}\}$ be independent standard Poisson processes. For each $i=1,\ldots,q$ and for each nonempty subset $J$ of $K_i$, let $\{e(i,J,l), l\in\mathbb{N}\}$ be a sequence of i.i.d. random variables distributed uniformly on $J$, independent across $i$ and $J$. Assume also that the four random objects $(N,\{\widetilde{\mu}_k,\widehat{\mu}_k\},X_0,\{B_{k,0}\})$, $A$, $\{S^{(i)}\}$ and $\{e(i,J,l)\}$ are mutually independent. Define*

(3.1) $$D^{(i)}(t) = S^{(i)}(T^{(i)}(t)), \qquad i=1,\ldots,q,$$

*where*

$$T^{(i)}(t) = \sum_{k \in K_i} T_k(t), \qquad i = 1, \ldots, q, \tag{3.2}$$

*and consider*

$$D_k(t) = \sum_{s \in (0,t]: \Delta D^{(i)}(s) = 1} 1_{\{e(i, \{p \in K_i : B_p(s-) = 1\}, D^{(i)}(s)) = k\}}, \tag{3.3}$$

$$k \in K_i, i = 1, \ldots, q$$

*as a substitute for equation* (2.7). *Then the process* $\Sigma_1'$, *defined analogously to* $\Sigma_1$, *with* (3.1)–(3.3) *in place of* (2.7), *is equal in law to* $\Sigma_1$. *Similarly, the process* $\Sigma_2'$ *defined analogously to* $\Sigma_2$, *with* (3.1)–(3.3) *in place of* (2.7), *is equal in law to* $\Sigma_2$.

This proposition is basically a statement about superposition of Poisson processes. We omit the elementary proof.

Let $S$ denote a standard Poisson process, and set

$$\widehat{A}(t) = n^{-1/2}(A(t) - \lambda^n t), \qquad t \geq 0, \tag{3.4}$$

$$\widehat{S}(t) = n^{-1/2}(S(nt) - nt), \qquad t \geq 0. \tag{3.5}$$

We will use in what follows the well-known fact that both $\widehat{A}$ and $\widehat{S}$ converge weakly to a zero mean Brownian motion with diffusion coefficient $\lambda^{1/2} C_{\check{U}}$, and respectively, 1 (see e.g. Lemmas 2 and 4(i) of [2]).

PROOF OF THEOREM 2.1. Let $\varepsilon > 0$ be given and let $q \in \mathbb{N}$ and $\mu^{(i)} \in \mathbb{R}_+$, $i = 1, \ldots, q$ be numbers satisfying the following conditions:

- $\mu^{(1)} = 0$, $\mu^{(q)} \geq 1$,
- $0 < \mu^{(i)} - \mu^{(i-1)} \leq \varepsilon$, $i = 1, \ldots, q$,
- $\int_{[\mu^{(q)}, \infty)} x^2 \, d\widetilde{m} \leq \varepsilon$,
- for $i = 2, \ldots, q$, $\mu^{(i)}$ is a continuity point of $x \mapsto \widetilde{m}([0, x]) \equiv P(\widetilde{\mu}_1 \leq x)$.

Set $\mu^{(q+1)} = \infty$, and denote $r_i = [\mu^{(i)}, \mu^{(i+1)})$ for $i = 1, \ldots, q$. Let

$$K_i = \{k \in \{1, \ldots, N\} : \mu_k \in r_i\}, \qquad i = 1, \ldots, q. \tag{3.6}$$

Using equations (2.4) and (2.9), and Proposition 3.1 with the above $K_i$, write

$$\widehat{X}(t) = \widehat{X}_0 + n^{-1/2} A(t) - n^{-1/2} \sum_{i=1}^{q} S^{(i)}(T^{(i)}(t)) \tag{3.7}$$

$$= \widehat{X}_0 + W(t) + bt + F(t), \tag{3.8}$$

where

$$W(t) = W^{n,\varepsilon}(t) = \widehat{A}(t) - \sum_{i=1}^{q} W^{(i)}(t), \tag{3.9}$$

$$W^{(i)}(t) = n^{-1/2}(S^{(i)}(T^{(i)}(t)) - T^{(i)}(t)), \qquad i = 1, \ldots, q, \tag{3.10}$$

$$b = b(n) = n^{-1/2}(\lambda^n - n\lambda) \tag{3.11}$$
$$- n^{-1/2} \sum_{k=1}^{N} (\widetilde{\mu}_k - \mu) - n^{-1} \sum_{k=1}^{N} \widehat{\mu}_k - \mu \widehat{N},$$

$$F(t) = F^{n,\varepsilon}(t) = n^{-1/2} \sum_{i=1}^{q} \left( \sum_{k \in K_i} \mu_k t - T^{(i)}(t) \right) \tag{3.12}$$
$$= n^{-1/2} \int_0^t \sum_{k=1}^{N} \mu_k I_k(s) \, ds.$$

Denote

$$I^{(i)}(t) = \sum_{k \in K_i} I_k(t), \qquad i = 1, \ldots, q. \tag{3.13}$$

Let $t$ be such that $H(t) > 0$. By (2.15) and (2.16), and since departures from different servers do not occur at the same time, we have that $I_k(t-) = 1$ if and only if $H_k(t) \in [H(t), t)$, which in turn occurs if and only if $\Delta D_k(s) = 1$ for some $s \in [H(t), t)$. Let $\widetilde{k}$ be such that $H(t) = H_{\widetilde{k}}(t)$. Note that, for a given $k$, it is impossible to have that $\Delta D_k(s) = 1$ for more than one $s \in [H(t), t)$, since this would mean that a departure from server $k$ and a routing to this server have taken place some time at $s_1$ and, respectively, $s_2$, where $H(t) < s_1 < s_2 < t$, thus by (2.18) and (2.19), $H_k(s_2) \le H_{\widetilde{k}}(s_2) = H(t)$, which contradicts $H_k(s_2) \ge s_1 > H(t)$. Writing $\#C$ for the cardinality of a set $C$, we therefore have for $t$ such that $H(t) > 0$,

$$I^{(i)}(t-) = \#\{k \in K_i : H_k(t) \in [H(t), t)\}$$
$$= \#\{k \in K_i : \Delta D_k(s) = 1 \text{ for some } s \in [H(t), t)\}$$
$$= D^{(i)}(t-) - D^{(i)}(H(t)-).$$

In case when $H(t) = 0$ one has to add those servers that were idle throughout the interval $[0, t)$. Thus, with the convention that, for any process $L$, $L(H(t)-) = L(0)$ if $H(t) = 0$, we have $I^{(i)}(t-) = D^{(i)}(t-) - D^{(i)}(H(t)-) + e_0^{(i)}(t)$, where

$$e_0^i(t) = \#\{k \in K_i : I_k(s) = 1 \text{ for all } s \in [0, t)\}.$$

Note that

$$(3.14) \quad \sum_{i=1}^{q} \widehat{e}_0^i(t) := \sum_{i=1}^{q} n^{-1/2} e_0^i(t) \leq n^{-1/2} I(0) 1_{\{H(t)=0\}} = \widehat{X}_0^- 1_{\{H(t)=0\}}.$$

Denote $\widehat{S}^{(i)}(t) = n^{-1/2}(S^{(i)}(nt) - nt)$, $\widehat{I}^{(i)} = n^{-1/2} I^{(i)}$, $\widehat{I}_k = n^{-1/2} I_k$, $\widehat{B}_k = n^{-1/2} B_k$. Also, recall $h(t) = t - H(t)$ and let $N^{(i)} = \#K_i$. Then

$$(3.15) \quad \begin{aligned} \widehat{I}^{(i)}(t-) &= n^{-1/2}(S^{(i)}(T^{(i)}(t-)) - S^{(i)}(T^{(i)}(H(t)-))) + \widehat{e}_0^i(t) \\ &= E^{(i)}(t) + n^{-1/2} \mu^{(i)} N^{(i)} h(t), \end{aligned}$$

where

$$(3.16) \quad \begin{aligned} E^{(i)}(t) &= [W^{(i)}(t-) - W^{(i)}(H(t)-)] \\ &\quad + \sum_{k \in K_i} (\mu_k - \mu^{(i)}) \int_{H(t)}^{t} \widehat{B}_k(s) \, ds \\ &\quad - \mu^{(i)} \sum_{k \in K_i} \int_{H(t)}^{t} \widehat{I}_k(s) \, ds + \widehat{e}_0^i(t). \end{aligned}$$

Using (3.15), one can express $h$ in terms of $\sum_i \widehat{I}^{(i)}$, $\sum_i E^{(i)}$ and $\sum_i \mu^{(i)} N^{(i)}$, and then use (3.15) again to obtain an expression for each $\widehat{I}^{(i)}$ that does not involve $h$. Along with (3.12), this yields the following expression for $F$:

$$(3.17) \quad F(t) = \gamma \int_0^t \widehat{I}(s) \, ds + e(t),$$

where $e(t) = e_1(t) + e_2(t) + e_3(t) + e_4(t)$, and

$$e_1(t) = \sum_{i=1}^{q} \sum_{k \in K_i} (\mu_k - \mu^{(i)}) \int_0^t \widehat{I}_k(s) \, ds,$$

$$e_2(t) = \sum_{i=1}^{q} \mu^{(i)} \int_0^t E^{(i)}(s) \, ds,$$

$$e_3(t) = -\frac{\sum_{i=1}^{q} (\mu^{(i)})^2 N^{(i)}}{\sum_{i=1}^{q} \mu^{(i)} N^{(i)}} \sum_{i=1}^{q} \int_0^t E^{(i)}(s) \, ds,$$

$$e_4(t) = \left( \frac{\sum_{i=1}^{q} (\mu^{(i)})^2 N^{(i)}}{\sum_{i=1}^{q} \mu^{(i)} N^{(i)}} - \gamma \right) \int_0^t \widehat{I}(s) \, ds,$$

with the convention that when the denominator in the expressions for $e_3$ and $e_4$ is zero, one lets $e_3(t) = 0$ and $e_4(t) = -\gamma \int_0^t \widehat{I}(s) \, ds$. Hence by (3.8) and (2.12), we have

$$(3.18) \quad \widehat{X}(t) = \widehat{X}_0 + W(t) + bt + \gamma \int_0^t \widehat{X}(s)^- \, ds + e(t).$$

Recall that, by construction, the partition (3.6) depends on $\varepsilon$, and as a result, so do the processes $\widehat{X}$ [given in equation (3.7)] and $W$ [defined in equation (3.9)]. However, Proposition 3.1 asserts that, for each $n$, the law of $\Sigma_1'$ does not depend on the partition, and since $W$ is defined only in terms of the processes $A$, $\sum_{i=1}^{q} S^{(i)}(T^{(i)}(\cdot)) = \sum_{k=1}^{N} D_k$ and $T = \sum_{i=1}^{q} T^{(i)} = \sum_{k=1}^{N} T_k = \sum_{k=1}^{N} \mu_k \int_0^{\cdot} B_k(s)\,ds$, we see that the joint law of $\widehat{X} = \widehat{X}^{n,\varepsilon}$ and $W = W^{n,\varepsilon}$ does not depend on $\varepsilon$. This observation is used several times below.

If $a(n,\varepsilon)$, $n \in \mathbb{N}$, $\varepsilon \in (0,1)$ is a family of real-valued random variables, we write $a \in \mathcal{N}$ if *for every positive $\delta$ and $\delta'$ there exists an $\varepsilon$ such that*

$$\limsup_{n \to \infty} P(|a(n,\varepsilon)| > \delta) < \delta'.$$

Fix $\bar{t} > 0$ throughout. Define

$$\theta = \theta(n,\varepsilon) = \inf\{t : |e(t)| \geq 1\} \wedge \bar{t}.$$

Our main estimate is the following:

LEMMA 3.1. (i) *One has $b \Rightarrow \beta$.*

(ii) *Fix $\varepsilon$. Then for $i = 1, \ldots, q$, $\sup\{|n^{-1}T^{(i)}(t) - \rho_i t| : t \leq \theta\} \to 0$ in probability, as $n \to \infty$, where $\rho_i = \int_{r_i} x\,d\widetilde{m}$.*

(iii) *For $j = 1, 2, 3, 4$, $|e_j|_{*,\theta} \in \mathcal{N}$.*

The proof of this lemma appears later. We can now complete the proof of the theorem. Note first that, by Lemma 3.1(iii), $|e|_{*,\theta} \in \mathcal{N}$. However, since the joint law of $b$, $X$ and $W$ does not depend on $\varepsilon$, it follows from (3.18) that neither does the law of $e$. As a result, the statement $|e|_{*,\theta} \in \mathcal{N}$ simply asserts that, for every $\varepsilon$, $|e|_{*,\theta} \to 0$ in probability. By definition of $\theta$, we therefore have that, for every $\varepsilon$,

$$(3.19) \qquad P(\theta < \bar{t}) \to 0 \qquad \text{as } n \to \infty,$$

and thus $|e|_{*,\bar{t}} \to 0$ in probability, as $n \to \infty$. In what follows, we fix $\varepsilon$.

By (3.10) and the definition of $\widehat{S}^{(i)}$,

$$(3.20) \qquad W^{(i)}(t) = \widehat{S}^{(i)}(n^{-1}T^{(i)}(t)), \qquad i = 1, \ldots, q.$$

From Lemma 3.1(ii) we now have that $n^{-1}(T^{(1)}, \ldots, T^{(q)}) \to \widetilde{\rho}$ in probability, uniformly on $[0,\bar{t}]$, where $\widetilde{\rho}(t) = (\rho_1 t, \ldots, \rho_q t)$. Recall that $(\widehat{A}, \widehat{S}^{(1)}, \ldots, \widehat{S}^{(q)})$ are mutually independent, and that $\widehat{S}^{(i)}$ (resp., $\widehat{A}$) converges to a standard Brownian motion (a zero-mean Brownian motion with diffusion coefficient $\lambda^{1/2}C_{\widetilde{U}}$). Thus (3.9), (3.20) and the lemma on random change of time [3, p. 151] show that $W$ converges weakly to $\sigma w$, in the uniform topology on $[0,\bar{t}]$, where $w$ is a standard Brownian motion and $\sigma^2 = \lambda C_{\widetilde{U}}^2 + \mu$.

Finally, by the Skorohod representation theorem, we can assume without loss of generality that the random variables $\widehat{X}_0$, $b$, $\xi_0$ and $\beta$, and the processes $W$ and $w$ are realized in such a way that

(3.21) $\qquad (\widehat{X}_0, b, W) \to (\xi_0, \beta, \sigma w) \qquad$ in probability, as $n \to \infty$.

Let $\xi$ be the unique strong solution to equation (2.20) with data $(\xi_0, \beta, w)$. Combining (2.20) and (3.18), the inequality $|x^- - y^-| \le |x - y|$, and Gronwall's inequality shows

$$|\widehat{X} - \xi|_{*,\bar{t}} \le (|\widehat{X}_0 - \xi_0| + |b - \beta| + |W - \sigma w|_{*,\bar{t}} + |e|_{*,\bar{t}}) \exp(\gamma \bar{t}).$$

By (3.21) and the uniform convergence of $e$ to zero, we have shown that $\widehat{X}$ converges to $\xi$ in probability, uniformly on $[0, \bar{t}]$. Since $\bar{t}$ is arbitrary, this shows $\widehat{X} \Rightarrow \xi$, and the result follows. $\square$

PROOF OF LEMMA 3.1. We begin by showing that $b \Rightarrow \beta$. By the assumption (2.3) on $N$, and the central limit theorem and the law of large numbers for i.i.d. random variables, the second and third terms in (3.11) converge to the random variable $-\zeta$ and, respectively, the constant $-\widehat{\mu}$ from the statement of the theorem. Thus by (2.5) and (2.3), $b \Rightarrow \beta$. This proves part (i) of the lemma.

Recall that $W^{(i)}$ are given in (3.20). We show that the random variables $\{|W^{(i)}|_{*,\bar{t}}, i = 1, \ldots, q, n \in \mathbb{N}\}$ are tight. By (2.1) and (2.8), for $i \in \mathbb{N}$,

$$
(3.22) \quad n^{-1} T^{(i)}(t) = n^{-1} \sum_{k \in K_i} \widetilde{\mu}_k t + n^{-3/2} \sum_{k \in K_i} \widehat{\mu}_k t \\
- n^{-1} \sum_{k \in K_i} \mu_k \int_0^t I_k(s) \, ds.
$$

Note that $N^{(i)}/n \Rightarrow \widetilde{m}(r_i)$ as $n \to \infty$. Since by construction each $\mu^{(i)}$ is a continuity point of $x \mapsto \widetilde{m}[0, x]$, using the law of large numbers we have $l_i^n := n^{-1} \sum_{k: \widetilde{\mu}_k \in r_i} \widetilde{\mu}_k \Rightarrow \int_{r_i} x \, d\widetilde{m} = \rho_i$. Denoting by $u_i^n t$ the first term on the r.h.s. of (3.22), it is easy to see by (2.1) that $l_i^n - u_i^n$ converge weakly to zero. For $t \le \bar{t}$, the second term is bounded by $n^{-1/2} \bar{\mu} \bar{t}$ for all $n$ sufficiently large. Since the third term is negative, we have from (3.20) that $|W^{(i)}|_{*,\bar{t}} \le |\widehat{S}^{(i)}|_{*,u_i^n+1}$. By the tightness of $\widehat{S}^{(i)}$ (as processes that converge to a Brownian motion) and the tightness of $u_i^n$, we have the tightness of $|W^{(i)}|_{*,\bar{t}}$.

Next, note that $|e|_{*,\theta} \le 1$. Thus the tightness of the random variables $\widehat{X}_0$, $b$, $|W^{(i)}|_{*,\bar{t}}$ and $|\widehat{A}|_{*,\bar{t}}$, $n \in \mathbb{N}$ (as follows from the convergence of $\widehat{A}$), and an application of Gronwall's lemma on (3.18), by which

$$|\widehat{X}|_{*,\theta} \le (|\widehat{X}_0| + |W|_{*,\theta} + |b|\bar{t} + 1) \exp(\gamma \bar{t}),$$

imply that $\{|\widehat{X}|_{*,\theta}, n \in \mathbb{N}\}$ are tight. Since by (2.12), $\widehat{I} = \widehat{X}^-$, we also have that $\{|\widehat{I}|_{*,\theta}, n \in \mathbb{N}\}$ are tight.

Equipped with the tightness of $|\widehat{I}|_{*,\theta}$, we can prove part (ii) of the lemma, and moreover, show that the stopped processes $\{W^{(i)}(\cdot \wedge \theta), i = 1, \ldots, q, n \in \mathbb{N}\}$ are $C$-tight. Fix $\varepsilon$. For $i = 1, \ldots, q$, we show that the supremum over $t \leq \theta$ of the last term in (3.22) (in absolute value) converges to zero in probability, as $n \to \infty$. Indeed, using Cauchy–Schwarz inequality for each $s$ and the fact that $I_k(s)$ takes values in $\{0, 1\}$, this term is bounded, for $t \leq \theta$, by

$$n^{-1} \int_0^t \left(\sum_{k=1}^N \mu_k^2\right)^{1/2} I(s)^{1/2} \, ds \leq \left(n^{-1} \sum_{k=1}^N \mu_k^2\right)^{1/2} n^{-1/4} \bar{t} |\widehat{I}|_{*,\theta}^{1/2}.$$

Note that $n^{-1} \sum_{k=1}^N \mu_k^2 \Rightarrow \int x^2 \, d\widetilde{m} < \infty$, and thus the asserted estimate on the last term of (3.22) follows from the tightness of $|\widehat{I}|_{*,\theta}$. Denote by $\bar{w}$ the modulus of continuity

$$\bar{w}_\tau(x, \delta) = \sup_{|s-t| \leq \delta; s, t \in [0, \tau]} |x(s) - x(t)|, \qquad \delta > 0,$$

for $x : [0, \tau] \to \mathbb{R}$. By the discussion following display (3.22), we now have that $\sup\{|n^{-1} T^{(i)}(t) - \rho_i t| : t \leq \theta\} \to 0$ in probability as $n \to \infty$, and part (ii) of the lemma follows. Consequently, the convergence of $\limsup_{n \to \infty} P(\bar{w}_\theta(W^{(i)}, \delta) > \delta')$ to zero as $\delta \to 0$, for arbitrary $\delta'$, follows from (3.20) and the $C$-tightness of $\widehat{S}^{(i)}$. In view of the tightness of the random variables $|W^{(i)}|_{*,\theta}$ established above, this shows that for $i = 1, \ldots, q$, $\{W^{(i)}(\cdot \wedge \theta), n \in \mathbb{N}\}$ are $C$-tight.

We now show that

(3.23) $\qquad\qquad\qquad \{n^{1/2} |h|_{*,\theta}, n \in \mathbb{N}\}$ are tight.

To this end, consider the event $\{|h|_{*,\theta} > a\}$ for a given $n$, and a given positive constant $a$. On this event, there exists $t \leq \theta$ such that at least one server is idle throughout the time interval $[H(t), t)$ [recall the definition of $H$, and that $h(t) = t - H(t)$]. The routing policy under consideration routes customers to servers at the order at which servers become available. Consequently, if a certain server is idle throughout a given interval $[t_1, t_2)$ then all arrivals between $t_1$ and $t_2$ are necessarily routed to servers that were already idle at time $t_1$. On the event indicated above, it follows that the number of arrivals between the times $H(t)$ and $t$ is less than the number of idle servers at time $H(t)$. Hence, on this event,

$$A(t) - A(t - a) - 1 \leq A(t-) - A(H(t)-) \leq |I|_{*,\theta}.$$

Substituting $a_0 n^{-1/2}$ for $a$ (where $a_0$ is a constant) and noting that $\lambda^n \geq \lambda n / 2$ for all $n$ sufficiently large, we thus have by (3.4)

$$P(n^{1/2} |h|_{*,\theta} > a_0) \leq P\left(\inf_{t \leq \theta} [\widehat{A}(t) - \widehat{A}((t - a_0 n^{-1/2}) \vee 0)] \leq -1\right)$$

$$+ P(|\widehat{I}|_{*,\theta} > a_0\lambda/4),$$

for all sufficiently large $n$. Since, as we mentioned earlier, $\widehat{A}$ converge to a Brownian motion, they are $C$-tight. Therefore the first term on the r.h.s. above tends to zero as $n \to \infty$, for any fixed $a_0$. Thus, in view of the tightness of $\{|\widehat{I}|_{*,\theta}, n \in \mathbb{N}\}$ which we have proved, (3.23) follows.

We now estimate $e_1$. Given a time $t \leq \theta$, a particular server $k$ has completed $D_k(t)$ services by time $t$, and therefore there have been at most $D_k + 1$ idle periods by that time, that is, at most $D_k + 1$ time intervals $[u, u') \subset [0, t)$ such that $I_k(s) = 1$ for $s \in [u, u')$, and $\inf_{[v,v')} I_k = 0$ for any $[v, v') \subset [0, t)$, $[u, u') \subsetneq [v, v')$. Since the length of each idle period is bounded by $|h|_{*,\theta}$ [by definition of $H(t)$ and $h(t)$], we have from (2.7) and (2.8),

$$\int_0^t I_k(s)\,ds \leq (D_k(t) + 1)|h|_{*,\theta} \leq (S_k(\mu_k t) + 1)|h|_{*,\theta}, \qquad k = 1,\ldots,N,$$

provided $t \leq \theta$. Hence

$$e_1^q(t) := \int_0^t \sum_{k \in K_q} \mu_k \widehat{I}_k(s)\,ds \leq n^{-1/2}|h|_{*,\theta} \sum_{k \in K_q} \mu_k(S_k(\mu_k t) + 1).$$

Denoting

$$M = n^{-1} \sum_{k \in K_q} \mu_k(S_k(\mu_k t) + 1),$$

we have $E[M] \leq 2\bar{t} E[(\mu_1)^2 1_{\{\mu_1 \in r_q\}}]$, and using (2.1) and the third bullet in the definition of $\{\mu^{(i)}\}$, this can be bounded by $c\varepsilon$, where $c$ does not depend on $n$ or $\varepsilon$. Hence given $\alpha > 0$ and $\delta > 0$,

$$P(|e_1^q|_{*,\theta} > \delta) \leq P(n^{1/2}|h|_{*,\theta} > \alpha) + P(M > \delta/\alpha),$$

and the tightness of $n^{1/2}|h|_{*,\theta}$ and the fact that their law does not depend on $\varepsilon$ imply that $|e_1^q|_{*,\theta} \in \mathcal{N}$. Now $e_1$ is bounded above by

$$e_1^q + \sum_{i<q} \sum_{k \in K_i} \varepsilon \int_0^{\cdot} \widehat{I}_k(s)\,ds \leq e_1^q + \varepsilon \int_0^{\cdot} \widehat{I}(s)\,ds,$$

and so the tightness of $|\widehat{I}|_{*,\theta}$ and the fact that their law does not depend on $\varepsilon$ imply that $|e_1|_{*,\theta} \in \mathcal{N}$.

Consider now the expression (3.16) for $E^{(i)}$. Note that $\sup\{|W^{(i)}(t-) - W^{(i)}(H(t)-)| : t \leq \theta\} \leq \bar{w}_\theta(W^{(i)}, |h|_{*,\theta})$. Thus the $C$-tightness of the processes $W^{(i)}$ stopped at $\theta$ and (3.23) imply that the first term on the r.h.s. of (3.16), stopped at $\theta$, converges uniformly on $[0, \bar{t}]$ to zero in probability as

$n \to \infty$. As a result, so does $e_2^1 := \sum_{i=1}^{q} \mu^{(i)}[W^{(i)}(t-) - W^{(i)}(H(t)-)]$. Next,

$$e_2^2(t) := \sum_{i=1}^{q} \sum_{k \in K_i} \mu^{(i)}(\mu_k - \mu^{(i)}) \int_{H(t)}^{t} \widehat{B}_k(s) \, ds$$

$$\leq (n^{1/2}|h|_{*,\theta}) \left[ \varepsilon \sum_{i=1}^{q-1} \mu^{(i)} \frac{N^{(i)}}{n} + \frac{1}{n} \sum_{k \in K_q} \mu_k^2 \right].$$

The expression in square brackets converges weakly to $\varepsilon \sum_{i<q} \mu^{(i)} \widetilde{m}(r_i) + \int_{r_q} x^2 \, d\widetilde{m} \leq \varepsilon \int_{[0,\infty)} x \, d\widetilde{m} + \varepsilon$. Since $n^{1/2}|h|_{*,\theta}$ are tight and their law does not depend on $\varepsilon$, we see that $|e_2^2|_{*,\theta} \in \mathcal{N}$. Moreover, with

$$e_2^3(t) := \sum_{i=1}^{q} (\mu^{(i)})^2 \sum_{k \in K_i} \int_{H(t)}^{t} \widehat{I}_k(s) \, ds \leq (\mu^{(q)})^2 |\widehat{I}|_{*,\theta} |h|_{*,\theta}, \qquad t \leq \theta,$$

(3.23) and the tightness of $|\widehat{I}|_{*,\theta}$ imply that $|e_2^3|_{*,\theta}$ converges to zero in probability as $n \to \infty$. Finally, by (3.14) and (3.23),

$$\sum_{i=1}^{q} \mu^{(i)} \int_{0}^{\theta} \widehat{e}_0^i(s) \, ds \leq \mu^{(q)} |\widehat{X}_0| |h|_{*,\theta} \to 0 \qquad \text{in probability, as } n \to \infty.$$

The above estimate, along with the results regarding $e_2^1$, $e_2^2$ and $e_2^3$ imply that $|e_2|_{*,\theta} \in \mathcal{N}$.

Next, since $N^{(i)}/n \Rightarrow \widetilde{m}(r_i)$, it is clear that

$$(3.24) \quad c_j(n) := n^{-1} \sum_{i=1}^{q} (\mu^{(i)})^j N^{(i)} \Rightarrow \sum_{i=1}^{q} (\mu^{(i)})^j \widetilde{m}(r_i) =: \bar{c}_j, \qquad j = 1, 2.$$

Since $|\bar{c}_1 - \int x \, d\widetilde{m}| \leq \sum_{i<q} \varepsilon \widetilde{m}(r_i) + \int_{r_q} x \, d\widetilde{m} \leq 2\varepsilon$, and $|\bar{c}_2 - \int x^2 \, d\widetilde{m}| \leq \sum_{i<q} 2\varepsilon \int_{r_i} x \, d\widetilde{m} + \int_{r_q} x^2 \, d\widetilde{m} \leq 2\varepsilon \mu + \varepsilon$, we see that

$$(3.25) \qquad\qquad\qquad c_2(n)/c_1(n) - \gamma \in \mathcal{N}.$$

Along with the tightness of $|\widehat{I}|_{*,\theta}$ and the fact that its law does not depend on $\varepsilon$, (3.25) shows that $|e_4|_{*,\theta} \in \mathcal{N}$. An argument similar to the one for $e_2$ shows that $|\sum_{i=1}^{q} \int_{0}^{\cdot} E^{(i)}(s) \, ds|_{*,\theta} \in \mathcal{N}$. Clearly, this and (3.25) imply that $|e_3|_{*,\theta} \in \mathcal{N}$. This completes the proof of part (iii) of the lemma. □

PROOF OF PROPOSITION 2.1. By the way the policy is defined, if a server $k$ finishes an idle period at time $u$ then the length of its idle period is equal to $h^n(u-)$. Moreover, since $n^{1/2}(t_n - s) \to \infty$, it follows from the representation to the departure process (3.1), and Lemma 3.1(ii), that $\mathbf{K}_{s,t_n}^n$ is nonempty with probability tending to 1. Hence, with probability tending to 1,

$$\inf\{h(u-): s \le u \le t_n\} \le \underline{\mathbf{I}}_{s,t_n} \le \bar{\mathbf{I}}_{s,t_n} \le \sup\{h(u-): s \le u \le t_n\}.$$

Denoting $\widehat{h}(t) = n^{1/2} h(t)$ and $E(t) = \sum_i E^{(i)}(t)$, we have by (3.15)

$$\widehat{I}(t-) = E(t) + \widehat{h}(t) c_1(n),$$

where $c_1(n)$ is as in (3.24). Since by the proof of Theorem 2.1 and Lemma 3.1, $\widehat{I} \Rightarrow \xi^-$, and $c_1(n) - \mu \in \mathcal{N}$, and since $\xi$ has continuous sample paths, it suffices to prove that $E \to 0$ in probability, uniformly over $[s, t]$, for every $0 < s < t$. By (3.23), for fixed $s > 0$, $P(H(s) = 0, \theta \ge s) \to 0$ as $n \to \infty$. Hence by (3.14), one has $\sup\{\widehat{e}_0^i(u): s \le u \le t\} \to 0$ in probability. In conjunction with the estimates on the terms $e_2^1$, $e_2^2$ and $e_2^3$ in the proof of Lemma 3.1, this establishes the convergence of $E$ alluded to above, hence the result. □

PROOF OF THEOREM 2.2. Let $\varepsilon > 0$ be given. Let $\mu^{(1)} = (\mu_{\min} - \varepsilon) \vee 0$, $\mu^{(2)} \in (\mu_{\min}, \mu_{\min} + \varepsilon)$ and $\mu^{(3)} = \mu_{\max} + \varepsilon$, where $\mu^{(2)}$ is a continuity point of $x \mapsto \widetilde{m}([0, x])$. Let $r_1 = [\mu^{(1)}, \mu^{(2)})$ and $r_2 = [\mu^{(2)}, \mu^{(3)})$, and set

$$K_i = \{k \in \{1, \ldots, N\}: \mu_k \in r_i\}, \qquad i = 1, 2.$$

Equations (3.7)–(3.12) are valid, with $q = 2$. We keep the notation $I^{(i)}$ for $i = 1, 2$ (3.13), and $\widehat{I}^{(i)} = n^{-1/2} I^{(i)}$. It will be shown below that, given $\bar{t} > 0$ and $\varepsilon > 0$,

(3.26) $\qquad |\widehat{I}^{(2)}|_{*,\bar{t}} \to 0 \qquad$ in probability, as $n \to \infty$.

Based on (3.26), the argument presented here for Theorem 2.2 follows steps that are similar to those of the proof of Theorem 2.1, and are, in fact, much simpler. First, by (3.8) and (3.12), one has

(3.27) $\qquad \widehat{X}(t) = \widehat{X}_0 + W(t) + bt + \mu_{\min} \int_0^t \widehat{X}(s)^- \, ds + e(t),$

where

(3.28) $\qquad e(t) = \sum_{k=1}^N (\mu_k - \mu_{\min}) \int_0^t \widehat{I}_k(s) \, ds.$

Denote $\theta = \inf\{t > 0: |e(t)| \ge 1\} \wedge \bar{t}$. The property $b \Rightarrow \beta$ is proved as in the proof of Lemma 3.1. Also, the tightness of the random variables $\{|W^{(i)}|_{*,\bar{t}}, n \in \mathbb{N}\}$, and as a result, that of $\{|\widehat{X}|_{*,\theta}, n \in \mathbb{N}\}$ and $\{|\widehat{I}|_{*,\theta}, n \in \mathbb{N}\}$ is argued exactly as in the proof of Lemma 3.1. The supremum over $t \le \theta$ of the absolute value of the last term in (3.22) converges to zero in probability, since $\mu_k$ are assumed to be bounded and $|\widehat{I}|_{*,\theta}$ are tight. In view of this, the argument following (3.22) shows that, for $i = 1, 2$, $\sup\{|n^{-1}T^{(i)}(t) - \rho_i t| : t \leq \theta\} \to 0$ in probability, as $n \to \infty$.

Now, by (3.28), $|e|_{*,\theta} \leq 2\varepsilon \bar{t}|\widehat{I}|_{*,\theta} + (\mu_{\max} + \varepsilon)\bar{t}|\widehat{I}^{(2)}|_{*,\bar{t}}$, and thus by (3.26), $|e|_{*,\theta} \in \mathcal{N}$. Since the law of $e$ and $\theta$ does not depend on $\varepsilon$, we have that $|e|_{*,\theta} \to 0$ in probability, as $n \to \infty$, for any fixed $\varepsilon$. Given these facts, the completion of the proof is carried out precisely as in the proof of Theorem 2.1.

We now show (3.26). Note first that the probability of the event $\eta_1 := \{I^{(2)}(0) = 0\}$ converges to one as $n \to \infty$. Indeed, since $\widetilde{m}(r_1) > 0$ and $\mu^{(2)}$ is a continuity point of $x \mapsto \widetilde{m}([0, x])$, we have that $N^{(1)}/n \Rightarrow \widetilde{m}(r_1) > 0$. Moreover, $I(0)/n \Rightarrow 0$ by (2.3) and (2.12). Hence $P(\{k \in K_1 : B_{k,0} = 1\} \neq \varnothing) \to 1$ as $n \to \infty$, and by the assumption (2.22) we have that $P(\eta_1) \to 1$ as $n \to \infty$.

Given $u > 0$, consider now the event $\eta := \{|I^{(2)}|_{*,\bar{t}} > 2un^{1/2}\}$. On the event $\eta \cap \eta_1$ one can find $0 \leq s < t \leq \bar{t}$ such that $I^{(2)}(y) > 0$ for $y \in [s, t]$, and $I^{(2)}(t) - I^{(2)}(s) > un^{1/2}$. Since the servers in $K_2$ all have greater rate than those in $K_1$, the routing policy assigns all arrivals within $[s, t]$ to $K_2$ servers. Hence by (2.6), (2.7), (2.9), (2.11), (3.1), we have

$$un^{1/2} < I^{(2)}(t) - I^{(2)}(s) = D^{(2)}(t) - D^{(2)}(s) - A(t) + A(s),$$

and therefore

$$u < \widehat{S}^{(2)}(n^{-1}T^{(2)}(t)) - \widehat{S}^{(2)}(n^{-1}T^{(2)}(s)) - \widehat{A}(t) + \widehat{A}(s)$$
$$+ \sum_{k \in K_2} \mu_k \int_s^t \widehat{B}_k(y)\,dy - \lambda n^{1/2}(t-s) - n^{-1/2}(\lambda^n - \lambda)(t-s).$$

We have by (2.8) and (3.2) that $n^{-1}T^{(2)}(t) \leq 3\mu_{\max}\bar{t} =: \tau$. Also, by (2.5), the last term above is bounded by $c(t-s)$ for some constant $c$ independent of $n$ and $\varepsilon$. Hence on the event $\eta \cap \eta_1$, with $\delta = t - s$,

(3.29) $\quad u < \bar{w}_\tau(\widehat{S}^{(2)}, 2\mu_{\max}\delta) + \bar{w}_{\bar{t}}(\widehat{A}, \delta) + n^{1/2}C(n, \varepsilon)\delta + c\delta,$

where $C(n, \varepsilon) = n^{-1}\sum_{k \in K_2}\mu_k - \lambda$. By (2.1), (2.4), (2.3) and the definition of $K_2$, $C(n, \varepsilon) \to -g(\varepsilon)$ in probability, as $n \to \infty$, where $g(\varepsilon) = \int_{r_1} x\,d\widetilde{m} > 0$. We obtain

$$P(|\widehat{I}^{(2)}|_{*,\bar{t}} > 2u) = P(\eta) \leq P_1(n, \varepsilon, u) + P_2(n, \varepsilon, u) + P(\eta_1^c),$$

where

$P_j(n, \varepsilon, u) = P(\text{there exists } \delta \in \Delta_j \text{ such that } (3.29) \text{ holds}), \qquad j = 1, 2,$

and $\Delta_1 = (0, n^{-1/4}]$, $\Delta_2 = (n^{-1/4}, \bar{t})$. Given $\varepsilon$ and $u$, we have $\lim_n P_j = 0$ for $j = 1, 2$, by $C$-tightness of $\widehat{S}^{(2)}$ and $\widehat{A}$, and the strict negativity of the weak limit of $C$. Since $\lim_n P(\eta_1^c) = 0$ and $u > 0$ is arbitrary, statement (3.26) follows. This completes the proof of the theorem. $\square$

**4. Policy P1 in deterministic environment.** The main difficulty dealt with in proving Theorem 2.1 is perhaps not in the randomness of the service rates but in establishing the reduction to a one-dimensional process. We would therefore like to present a version of Theorem 2.1 in a deterministic environment that follows upon minor modifications. We consider here the following setting. We assume that $N^n = n$ for all $n$, and, in place of (2.1), model the service rates as deterministic, nonnegative constants

$$(4.1) \qquad \mu_k^n = \widetilde{\mu}_k^n + n^{-1/2}\widehat{\mu}_k^n, \qquad k = 1, \ldots, n.$$

The definition of the constants $(\widetilde{\mu}_k^n, \widehat{\mu}_k^n)$ is based on a given probability measure $m$, supported on a finite set $M \equiv \{M_1, \ldots, M_L\} \subset \mathbb{R}_+ \times \mathbb{R}$. Denote $M_l = (\widetilde{M}_l, \widehat{M}_l)$ for $l = 1, \ldots, L$. With $p_l := \widetilde{m}(M_l)$, $l = 1, \ldots, L$, the constants $(\widetilde{\mu}_k^n, \widehat{\mu}_k^n)$ are assumed to have values in $M$ and to satisfy

$$(4.2) \qquad |\#\{k : (\widetilde{\mu}_k^n, \widehat{\mu}_k^n) = M_l\} - np_l| \leq c, \qquad l = 1, \ldots, L, n \in \mathbb{N},$$

where $c$ is a constant that does not depend on $n$ and $k$. The marginals of $m$ are denoted by $\widetilde{m}$ and $\widehat{m}$, and, as in (2.2), $\widetilde{m}$ is assumed to satisfy $\mu := \int x\, d\widetilde{m} = \sum_l p_l \widetilde{M}_l \in (0, \infty)$. With the exception that the constants (4.1) replace the random variables (2.1), we keep here the complete setting of Section 1, including the definition of policy P1, that in the current setting will be denoted by P1′. Note that due to our assumption on $N^n$, the constant $\nu$ is zero.

COROLLARY 4.1. *Under policy P1′, the processes $\widehat{X}^n$ converge weakly to the unique solution to the equation*

$$\xi(t) = \xi_0 + \sigma w(t) + \beta t + \gamma \int_0^t \xi(s)^- ds, \qquad t \geq 0,$$

*where $\sigma^2 = \lambda C_{\check{U}}^2 + \mu$, $\beta = \widehat{\lambda} - \widehat{\mu}$, $\widehat{\mu} = \int x\, d\widehat{m} \equiv \sum_l p_l \widehat{M}_l$, $\gamma = \int x^2\, d\widetilde{m} / \int x\, d\widetilde{m} \equiv \sum_l p_l \widetilde{M}_l^2 / \sum_l p_l \widetilde{M}_l$, and $w$ is a standard Brownian motion, independent of $\xi_0$.*

PROOF. The proof proceeds as that of Theorem 2.1, with trivial modifications, given the observations below.

In Proposition 3.1, "a partition measurable on $\sigma\{N, \{\widetilde{\mu}_k, \widehat{\mu}_k\}\}$" reduces to "a deterministic partition."

With the notation of Section 3, we have by (4.2) that $n^{-1}N^{(i)} = \widetilde{m}(r_i) + O(n^{-1})$ as $n \to \infty$. Hence provided that $\varepsilon$ is sufficiently small we have $n^{-1}N^{(i)} \to \widetilde{m}(r_i)$, and for $j = 1, 2$, $n^{-1}\sum_{k:\widetilde{\mu}_k^n \in r_i}(\widetilde{\mu}_k^n)^j \to \int_{r_i} x^j\, d\widetilde{m}$, and moreover $n^{-1/2}\sum_{k=1}^n(\widetilde{\mu}_k^n - \mu) \to 0$, as $n \to \infty$. Similarly, $n^{-1}\sum_{k=1}^n \widehat{\mu}_k^n \to \int x\, d\widehat{m}$ as $n \to \infty$.

The third bullet in the definition of the partition (see the first paragraph in the proof of Theorem 2.1) forces that $\int_{[\mu^{(q)},\infty)} x^2\, d\widetilde{m} = 0$. This trivializes the treatment of $e_1^q$ in Lemma 3.1, since this term is zero. □

**Acknowledgment.** I would like to thank the referees for various useful comments, and for pointing out Remark 2.3.

DEPARTMENT OF ELECTRICAL ENGINEERING
TECHNION—ISRAEL INSTITUTE OF TECHNOLOGY
HAIFA 32000
ISRAEL
E-MAIL: atar@ee.technion.ac.il